\def\Max{\bigcirc\hspace{-4.54mm}\vee \,\,}

\def\Exp{\rm Exp}
\documentclass[12pt,a4paper,eqno]{amsart}
\usepackage{graphicx,epsfig}
\usepackage{amssymb}
\usepackage{mathtools}
\mathtoolsset{showonlyrefs}
\newenvironment{namelist}[1]{%
\begin{list}{}
     {
      
      \settowidth{\labelwidth}{#1}
      \setlength{\leftmargin}{1.1\labelwidth}
               }
      }{%
\end{list}}
\newtheorem{thm}{Theorem}[section]
\newtheorem{lem}[thm]{Lemma}
\newtheorem{prop}[thm]{Proposition}

\newtheorem{defn}[thm]{Definition}
\newtheorem{rem}[thm]{Remark}
\renewcommand{\uppercasenonmath}[1]{} 
\begin{document}
\title[Generalized Poisson processes]{\Large Classical definitions of the Poisson process do not coincide in the case of generalized convolutions}
\author[B.H. Jasiulis-Go{\l}dyn, J.K. Misiewicz]{\large B.H. Jasiulis-Go{\l}dyn$^1$ and J.K. Misiewicz $^2$}
\thanks{$^1$ Institute of Mathematics, University of Wroc{\l}aw, pl. Grunwaldzki 2/4, 50-384 Wroc{\l}aw, Poland, e-mail: jasiulis@math.uni.wroc.pl \\
$^2$ Faculty of Mathematics and Information Science, Warsaw
University of Technology, ul. Koszykowa 75, 00-662 Warsaw, Poland, e-mail:
j.misiewicz@mini.pw.edu.pl 
}
\maketitle
\vspace{-2mm}
\begin{abstract}
In the paper we consider a generalizations of the notion of Poisson process to the case when classical convolution is replaced by generalized convolution in the sense of K. Urbanik \cite{Urbanik64} following two classical definitions of Poisson process. First, for every generalized convolution $\diamond$ we define $\diamond$-generalized Poisson process type I as a Markov process with the $\diamond$-generalized Poisson distribution. Such processes have stationary independent increments in the sense of generalized convolution, but usually they do not live on $\mathbb{N}_0$. The $\diamond$-generalized Poisson process type II is defined as a renewal process based on the sequence $S_n$, which is a Markov process with the step with the lack of memory property. Such processes take values in $\mathbb{N}_0$, however they do not have to be Markov processes, do not have to have independent increments, even in generalized convolution sense. It turns out that the second construction is possible only for monotonic generalized convolutions which admit the existence of distributions with lack of memory, thus we also study these properties.
\end{abstract}

\tableofcontents

\section{Introduction}

The paper deals with two methods of generalizing of the notion of Poisson process to the case when the classical convolution is replaced by a generalized convolution in the sense of K. Urbanik \cite{Urbanik64}. As far as we know the generalized Poisson processes in the sense of generalized convolution have not been studied yet, but Markov processes which can be considered as additive processes in the sense of generalized convolutions have been considered in \cite{WeakLevyProc, Thu94, Thu}. Properties of $\diamond$-generalized Poisson distribution were studied in \cite{reprLevy, Urbanik64}.  Basic properties and a list of examples of generalized convolutions is given in Section 2. Then we consider two possible definitions of Poisson process which in the case of classical convolution coincides:

{\bf I.} The first definition states that the stochastic process $\{N_t\colon t \geqslant 0\}$ is a Poisson process if it starts from zero ($N_0 = 0$ a.e.), it has independent and stationary increments and the distribution of $N_t$ is the Poisson distribution with parameter $a t$, $a >0$, i.e.
$$
N_t \sim \exp \left( a t \delta_1 \right).
$$
In Section 3, replacing classical convolution by a generalized one $\diamond$ and  the measure $\exp \left( a t \delta_1 \right)$ by its equivalent  $\Exp_{\diamond} \left( a t \delta_1 \right)$ we obtain a Markov process with transition probabilities $\delta_x \diamond \Exp_{\diamond} \left( a (t-s) \delta_1 \right)$ which is an additive process with respect to convolution $\diamond$. In Section 3 we show that for  each generalized convolution $\diamond$ this construction is possible, the obtained processes have stationary and independent increments in the sense of convolution $\diamond$, however most of them do not take values in natural numbers. We calculated explicitly the infinitesimal operators for such processes.

{\bf II.} The second classical definition of the Poisson process is based on direct construction as a renewal process via the sequence $\{T_n \colon n \in \mathbb{N} \}$ of i.i.d. random variables with the exponential distribution $\Gamma(1, a)$, $a>0$. Then, for $S_0=0$ and $S_n = T_1 + \dots + T_n$ for $n \geqslant 1$, we have
$$
N_t = \left\{ \begin{array}{l}
 \inf\{ n \colon S_{n+1} > t \} \\[2mm]
 \infty \; \hbox{ if such $n$ does not exists }.
 \end{array} \right.
$$
Following this construction we shall first construct an increasing sequence of change times $S_n$ as a Markov process with the step distribution $\nu$ having lack of memory property with respect to convolution $\diamond$ and the transition probability $\delta_x \diamond \nu$. It turned out the only some of generalized convolutions are monotonic, so that the sequence $S_n$ is monotonically increasing - description is given in Section 4. Then, in Section 5, we show that existence of the distribution with the lack of memory property for given convolution $\diamond$ depends on the regularity of convolution. Finally in Section 6 we construct $\diamond$-generalized Poisson process as a renewal process for generalized convolutions which are monotonic and admit distributions with lack of memory. On two examples for stable convolution and for Kendall convolution we show that such processes do not have to be Markov processes, do not have to have independent increments, even with respect to generalized convolution. However such processes are taking values  in natural numbers.

\section{Urbanik's generalized convolutions}
The generalized convolutions on the set $\mathcal{P}_+$ of probability measures on the Borel subsets of the positive half line were defined by Urbanik (see \cite{Urbanik64}). It is also possible to consider more general binary operation also called generalized convolution defined on the set $\mathcal{P}$ of all probability measures on the Borel subsets of the real line $\mathbb{R}$ (see \cite{Basia,conv1,MOU}) or generalized convolutions on the set $\mathcal{P}_s$ of all symmetric probability measures on $\mathbb{R}$. When our result holds in every of these cases we say that considered  measures are living on $\mathbb{K}$.

For simplicity we will use notation $T_a$ for the rescaling operator (dilatation of probability measure) defined by $(T_a\lambda)(A) = \lambda ({A/a})$ for every Borel set $A$ when $a \neq 0$, and $T_0 \lambda = \delta_0$.

\begin{defn}\label{df:1}
A commutative and associative $\mathcal{P}$-valued binary
operation $\diamond$ defined on $\mathcal{P}_{+}^2$ is called a
\emph{generalized convolution} if for all $\lambda,\lambda_1,\lambda_2 \in \mathcal{P}_{+}$ and $a \geqslant 0$ we have:
\begin{namelist}{llll}
\item[{\rm (i)}] $\delta_0 \diamond \lambda = \lambda$ ; %
\item[{\rm (ii)}] $(p\lambda_1 +(1-p)\lambda_2) \diamond \lambda = p ( \lambda_1 \diamond \lambda) + (1-p)(\lambda_2
          \diamond \lambda) $ whenever $p\in [0,1]$; %
\item[{\rm (iii)}] $T_a(\lambda_1 \diamond \lambda_2) = (T_a \lambda_1) \diamond (T_a \lambda_2)$ ; %
\item[{\rm (iv)}] if $\lambda_n \rightarrow \lambda$ then $\lambda_n \diamond \eta \rightarrow \lambda \diamond
           \eta$  for all $\eta \in \mathcal{P}$ and $\lambda_n \in \mathcal{P}_{+}$ ,
\end{namelist}
where $\rightarrow$ denotes weak convergence of probability measures.
\end{defn}
In the Urbanik definition (see \cite{Urbanik64}) of generalized convolutions on $\mathcal{P}_+$ there is one more condition connected with Central Limit Theorem:
\begin{namelist}{llll}
\item[{\rm (v)}] there exists a sequence $(c_n)_{n\in\mathbb{N}}$ of positive numbers such that the sequence $T_{c_n}
           \delta_1^{\diamond n}$ converges to a measure different from $\delta_0$.
\end{namelist}
The wide discussion on condition (v) for generalized convolution on $\mathcal{P}$ one can find in \cite{Basia}. In particular it was proven there that there exists generalized convolution without property (v).

The set $(\mathcal{P}_{+}, \diamond)$ we call a \emph{generalized convolution algebra}. A continuous mapping $h: \mathcal{P} \to \mathbb{R}$ such that
\begin{itemize}
\item $h(p\lambda + (1-p)\nu)= p h(\lambda) + (1-p) h(\nu)$,
\item $h(\lambda  \diamond \nu)=h(\lambda) h(\nu)$
\end{itemize}
for all $\lambda, \nu \in \mathcal{P}_{+}$ and $p \in (0,1)$, is called a \emph{homomorphism of $(\mathcal{P}_{+}, \diamond)$}. Every convolution algebra $(\mathcal{P}_{+}, \diamond)$ admits two trivial homomorphisms: $h \equiv 1$ and $h \equiv 0$. We say that a generalized convolution is \emph{regular} if it admits non-trivial homomorphism.  If the generalized convolution is regular then its homomorphism is uniquely determined in the sense that if $h_1, h_2$ are homomorphisms of $(\mathcal{P}_{+}, \diamond)$ then there exists $c>0$ such that $h_1(\lambda) = h_2(T_c \lambda)$ (for details see \cite{Urbanik64}). For every probability measure $\lambda \in \mathcal{P}_+$ we have
$$
h(\lambda) = \int_0^{\infty} h(\delta_x) \lambda(dx).
$$
It was shown in \cite{Urbanik64} that the generalized convolution is regular if and only if there exists unique up to a scale function
$$
\mathcal{P}_{+} \ni \lambda \longrightarrow \Phi_{\lambda} \in C([0,\infty))
$$
such that $\Phi_{\lambda}(t) = 1$, and for all $\lambda, \nu, \lambda_n \in \mathcal{P}_{+}$ the following conditions hold:
\begin{namelist}{ll}
\item{1.} $\Phi_{p \lambda + q \nu}(t) = p \Phi_{\lambda}(t) + q \Phi_{\nu}(t)$, for $p,q \geqslant 0$, $p+q = 1$;
\item{2.} $\Phi_{\lambda \diamond \nu}(t) = \Phi_{\lambda}(t) \Phi_{\nu}(t)$;
\item{3.} $\Phi_{T_a \lambda}(t) = \Phi_{\lambda}(at)$ for $a \geqslant 0$;
\item{4.} the uniform convergence of $\Phi_{\lambda_n}$ on every bounded interval is equivalent to the weak convergence of $\lambda_n$
\end{namelist}
The function $\Phi_{\lambda}$ is called the $\diamond$-generalized characteristic function of the measure $\lambda$. Moreover, it was also shown in \cite{Urbanik64} that up to a scale parameter
$$
\Phi_{\lambda}(t) = \int_0^{\infty} h(T_t \delta_x) \lambda(dx).
$$
Since for all $\lambda_1,\lambda_2\in\mathcal{P}_+$ we have
    $$
    \lambda_1 \diamond \lambda_2(A)=\int_0^{\infty} \int_0^{\infty} \left( \delta_x \diamond \delta_y \right) (A)\, \lambda_1(dx) \lambda_2(dy),
    $$
every generalized convolution is uniquely determined by the {\em probability kernel}
$$
\rho_{x,y} := \delta_x \diamond \delta_y.
$$
Evidently for $x,y,c \geqslant 0$
\begin{itemize}
    \item
    $ \rho_{x,0}=\delta_x$, \quad \quad \quad  \quad $\bullet \, \rho_{x,y}=\rho_{y,x}$,
    \item
    $ T_c\rho_{x,y}=\rho_{cx,cy}$, \quad \quad
    $\bullet \, \rho_{x,y}=T_v\rho_{z,1}$, where $v=x\vee y,\; z=\frac{x\wedge y}{x\vee y}$.
        \end{itemize}

The origin of the generalized convolution we can also find in the Kingman paper (\cite{King}), where the first example of random walk under generalized convolution (called the Kingman or Bessel convolution) was considered. The Kingman convolution has the natural interpretation at the interference phenomena (see \cite{plastimo}). One can find many open problems connected with generalized convolutions (e.g. in \cite{vol2}). All the basic information about the examples listed bellow one can find in  \cite{Bingham, JasKula, King, KU2, Urbanik86, Urbanik84, Urbanik76, Urbanik73, Urbanik64, vol3}

\vspace{2mm}

\noindent
{\bf Examples.}
\begin{namelist}{ll}
\item[{\bf 0.}] The classical convolution is evidently an example of generalized convolution. It will be denoted simply by $\ast$:
    $$
    \delta_a \ast \delta_b = \delta_{a+b}, \quad \quad h(\delta_t)=e^{-t}.
    $$
\item[{\bf 1.}] The so-called {\em symmetric} generalized convolution on $\mathcal{P}_{+}$ is defined by
    $$
    \delta_a \ast_{s} \delta_b = \frac{1}{2}\, \delta_{|a-b|} + \frac{1}{2} \,\delta_{a+b}.
    $$
    The corresponding homomorphism is defined by $h(\delta_t) = e^{-t}$ or, if we consider symmetric distributions, by $h(\delta_t) = \cos{t}$.
    The name {\em symmetric} comes from the fact that this convolution can be easily extended to a generalized convolution on $\mathcal{P}^2$ taking values in the set of symmetric measures $\mathcal{P}_{s}$:
    $$
    \delta_a \ast_{s} \delta_b = \frac{1}{4} \,\delta_{a-b} + \frac{1}{4} \,\delta_{-a+b} + \frac{1}{4}\, \delta_{-a-b} \frac{1}{4} \, \delta_{a+b}.
    $$
\item[{\bf 2.}] Another generalized convolution (called by Urbanik $(\alpha,1)$-convolution, see \cite{Urbanik76}) we obtain in a similar way for every $\alpha>0$ defining
    $$
    \delta_a \ast_{s, \alpha} \delta_b = \frac{1}{2} \, \delta_{|a^{\alpha}-b^{\alpha}|^{1/{\alpha}}} + \frac{1}{2}\, \delta_{(a^{\alpha}+b^{\alpha})^{1/{\alpha}}}.
    $$
\item[{\bf 3.}] For every $\alpha \in (0,\infty]$ the formula
$$
\delta_a \ast_{\alpha} \delta_b = \delta_c, \quad \quad a,b \geqslant 0, \;\; c = \|(a,b)\|_{\alpha} = (a^{\alpha} + b^{\alpha})^{1/{\alpha}}
$$
defines a generalized convolution $\ast_{\alpha}$ (called $\alpha$-stable convolution) on $\mathcal{P}_{+}^2$ with the corresponding homomorphism $h(\delta_t) = e^{-t^{\alpha}}$.
\item[{\bf 4.}] The Kendall convolution $\vartriangle_{\alpha}$ on $\mathcal{P}_{+}^2$, $\alpha > 0$, is defined by
$$
\delta_x \vartriangle_{\alpha} \delta_1 = x^{\alpha} \pi_{2\alpha} + (1-x^{\alpha}) \delta_1, \quad x\in [0,1],
$$
where $\pi_{2\alpha}$ is the Pareto measure with the density function $\pi_{2\alpha}(x) = 2\alpha x^{-2\alpha -1} \mathbf{1}_{[1,\infty)}(x)$. The corresponding homomorphism given by $h(\delta_t) = (1 - t^{\alpha})_{+}$ is also the Archimedean copula generator and it is strictly connected with the Williamson transform (see \cite{Neslehova}). The characterization of the Kendall convolution one can find in \cite{MisJas2}.
\item[{\bf 5.}] The Kingman convolution $\otimes_{\omega_s}$ on $\mathcal{P}_{+}^2$, $s> - \frac{1}{2}$, is defined by
$$
\delta_a \otimes_{\omega_s} \delta_b = \mathcal{L}\left(\sqrt{ a^2 + b^2 + 2ab \theta_s }\right),
$$
where $\mathcal{L}(X)$ is the distribution of the random element $X$ and $\theta_s$ is a random variable with the density function
$$
f_s (x) = \frac{\Gamma(s+1)}{\sqrt{\pi}\, \Gamma(s + \frac{1}{2})} \bigl( 1 - x^2\bigr)_{+}^{s - \frac{1}{2}}.
$$
If $n:= 2(s+1) \in \mathbb{N}$, $n >1$, the variable $\theta_s$ can be interpreted as one dimensional projection $\omega_{n,1}$ of the uniform distribution $\omega_n$ on the unit sphere $S_{n-1} \subset \mathbb{R}^n$. If $n=1$ and $s = - \frac{1}{2}$ then $\theta_s$ has the discrete distribution $\frac{1}{2} \delta_{-1} + \frac{1}{2} \delta_1$. The homomorphism for the Kingman convolution is given by
 $$
    h(\delta_t) =  \Gamma\left( \frac{\beta}{2} \right) \left( \frac{2}{t} \right)^{\frac{\beta}{2} -1} J_{\frac{\beta}{2} -1} (t),
    $$
    where $J_{r}$ is the Bessel function of the first kind.
\item[{\bf 6.}] Max-convolution is defined by
$$
\delta_a \Max \delta_b = \delta_{\max\{a,b\}}
$$
and its homomorphism is $h(\delta_t) = \mathbf{1}_{[0.1]}(t)$. This generalized convolution is not regular. More about Max-convolution we can find in \cite{KU2,Urbanik73,  Urbanik76, Urbanik84, Urbanik86, Urbanik87}.
\item[{\bf 7.}] A combination of Kingman convolution and $(\alpha,1)$ convolution, called by Urbanik $(\alpha, \beta)$-convolution, for $0<\alpha <\infty, 0 < \beta < \infty$, is defined as
    $$
    \delta_a \otimes_{\alpha,\beta}\delta_b = \mathcal{L}\left( \bigl(a^{2\alpha} + b^{2\alpha} + 2 x^{\alpha} b^{\alpha} \theta \bigr)^{1/{2\alpha}} \right),
    $$
    where $\mathcal{L}(X)$ is the distribution of the random variable $X$ and $\theta$ is a random variable with the density function
    $$
    f_{(\beta-2)/2} (x) = \frac{\Gamma({\beta}/2)}{\sqrt{\pi}\, \Gamma({(\beta-1)}/2)} \left( 1 - x^2 \right)_{+}^{(\beta-3)/2}.
    $$
The homomorphism is given by
 $$
    h(\delta_t) =  \Gamma\left( \frac{\beta}{2} \right) \left( \frac{2}{t^{\alpha}} \right)^{\frac{\beta}{2} -1} J_{\frac{\beta}{2} -1} (t^{\alpha}).
    $$
\item[{\bf 8.}] The Kucharczak convolution $\bigcirc \hspace{-3.2mm} \alpha$, $\alpha \in (0,1)$, is defined by
    $$
    \delta_a {\bigcirc \hspace{-3.2mm} \alpha}\,\, \delta_b  (dx) = \frac{a^{\alpha} b^{\alpha} \sin(\pi \alpha)(2x - a -b)}{\pi\,(x-a-b)^{\alpha} (x-a)^{\alpha} (x-b)^{\alpha}} \, \mathbf{1}_{[(a^{\alpha} + b^{\alpha})^{1/{\alpha}}, \infty)}(x) dx.
    $$
    In this case
    $$
    h(\delta_t) = \Gamma(\alpha)^{-1} \Gamma(\alpha, t),
    $$
    where $\Gamma(\alpha, t)$ is the incomplete Gamma function.
\item[{\bf 9.}] The Vol'kovich convolution $\vartriangle_{1,\beta}$ for $0<\beta< \frac{1}{2}$ (see \cite{vol3,vol1}) is given by
    $$
    \delta_a \vartriangle_{1,\beta} \delta_b (dx) = \frac{2 a^{2\beta} b^{2\beta}dx}{B(\beta, \frac{1}{2} - \beta)} \left( \bigl( x^2 - (a-b)^2 \bigr)_{+} \bigl( (a+b)^2 - x^2 \bigr)_{+} \right)^{-\beta - \frac{1}{2}}.
    $$
    In this case
    $$
    h(\delta_t)  = \frac{2^{1-\beta} t^{\beta}}{\Gamma(\beta)}\, K_{\beta}(t),
    $$
    where $K_{\beta}$ is the MacDonald function given by
    $$
    K_{\beta}(t) = \frac{\sqrt{\pi}}{\Gamma(\beta + \frac{1}{2})} \left(\frac{t}{2}\right)^{\beta} \int_0^{\infty} e^{-t {\rm ch}s} {\rm sh}^{2\beta}(s) ds.
    $$
\item[{\bf 10.}] In \cite{KU2} for $\alpha \in (0,1)$ authors considered the following measure:
    $$
    \mu = \left(2 - 2^{-\alpha}\right) \sum_{n=0}^{\infty} 2^{-1-n(\alpha +1)} T_{2^n} \pi_{\alpha},
    $$
    where $\pi_{\alpha}$ is the Pareto distribution with the density function $\pi_{\alpha}(x) = \alpha x^{-\alpha - 1} \mathbf{1}_{[1, \infty)}$.
    They proved that for every pair $a,b>0$ there exists a unique probability measure $\varrho(a,b)$ fulfilling the equality
    $$
    T_a \mu \Max T_b \mu = \mu \circ \varrho(a,b).
    $$
    Setting $\delta_a \nabla_{\alpha} \delta_b := \varrho(a,b)$ they've got a generalized convolution. In a similar way many other generalized convolutions can be constructed on a base of known convolutions (see e.g. \cite{JasKula}). The corresponding homomorphism is given by
    $$
    h(\delta_t) = \left( 1 - 2^{(1+\alpha)[\log_2t]} - (2 - 2^{-\alpha} ) (1- 2^{[\log_2t]}) t^{\alpha} \right) \mathbf{1}_{[0,1]}(t),
    $$
    where the square brackets denote the integer part.
\end{namelist}

\section{$\diamond$-Generalized Poisson process of type I}
Following the classical construction of the Poisson process given by definition I we introduce the generalized Poisson process of the first kind based on the generalized exponent of the Dirac-delta distribution, where
\begin{defn}
For $\lambda \in \mathcal{P}$ and $a>0$ the $\diamond$-generalized compound Poisson distribution ${\rm Exp}_{\diamond} ( a \lambda)$ with respect to the generalized convolution $\diamond$ is given by
$$
{\rm Exp}_{\diamond} ( a \lambda) = e^{-a} \sum_{k=0}^{\infty} \frac{a^k}{k!} \lambda^{\diamond k}.
$$
\end{defn}
The $\diamond$-characteristic function for the measure ${\rm Exp}_{\diamond} ( a \lambda)$ for regular generalized convolution can be calculated as follows:
\begin{eqnarray*}
\lefteqn{\int_0^{\infty} h(T_t \delta_x) {\rm Exp}_{\diamond} ( a \lambda)(dx) = e^{-a} \sum_{k=0}^{\infty} \frac{a^k}{k!}  \int_0^{\infty} h(T_t \delta_x) \lambda^{\diamond k}(dx) } \\
 & = & e^{-a} \sum_{k=0}^{\infty} \frac{a^k}{k!} \left[ \int_0^{\infty} h(T_t \delta_x) \lambda(dx)\right]^k = \exp \left\{ - a \left( 1 - \Phi_{\lambda}(t)\right) \right\}.
\end{eqnarray*}
Using this expression it is easy to see that
$$
{\rm Exp}_{\diamond} ( a \lambda) \diamond {\rm Exp}_{\diamond} ( b \lambda) = {\rm Exp}_{\diamond} (( a+b) \lambda),
$$
which shows that the measure ${\rm Exp}_{\diamond} ( a \lambda)$ is infinitely divisible in the sense of generalized convolution $\diamond$. More about $\diamond$-infinite divisibility of measures (called sometimes also infinite decomposability) one can find in \cite{Urbanik64, reprLevy}.

\begin{defn}
If $\diamond$ is a regular generalized convolution then by the generalized Bernoulli distribution (notation $\mathcal{GB}(\diamond, n, p)$) we understand the following measure
$$
\nu_n := \left( p \delta_1 + (1-p) \delta_0 \right)^{\diamond n} = \sum_{k=0}^n {{n}\choose{k}} p^k (1-p )^{n-k} \delta_1^{\diamond k}.
$$
\end{defn}
The next proposition shows that, similarly to the classical convolution, it is a weak limit of the generalized Bernoulli distributions.
\begin{prop}
Let  $\nu_n = \mathcal{GB}(\diamond, n, p_n)$ be such that $n p_n \rightarrow a>0$ for $n \rightarrow \infty$ and let $(\mathcal{P}_{+}, \diamond)$ be regular. Then
$$
\nu_n \rightarrow {\rm Exp}_{\diamond} ( a \delta_1), \quad n\rightarrow \infty,
$$
where $\rightarrow$ denotes weak convergence of distributions.
\end{prop}
\noindent
{\bf Proof.} Calculating the generalized characteristic function for the measure $\nu_n$ we have
\begin{eqnarray*}
h \left( T_t \nu_n\right) & = & \left( (1-p_n) h(\delta_0) + p_n \, h(\delta_t) \right)^n = \left( 1  - \frac{n p_n \,(1- h(\delta_t))}{n} \right)^n \\
 && \stackrel{n\rightarrow\infty}{\longrightarrow} \exp \left\{ - a (1- h(\delta_t))\right\} = h(T_t {\rm Exp}_{\diamond}(a \delta_1)),
\end{eqnarray*}
where $h$ is the nontrivial homomorphism for $\diamond$.
This means that $\nu_n \rightarrow {\rm Exp}_{\diamond}(\lambda \delta_1)$ by Th. 6 in \cite{Urbanik64} \qed \\

\vspace{2mm}

{\bf Examples.}
\begin{namelist}{ll}
\item[{\bf 1a.}] For the symmetric convolution we have
\begin{eqnarray*}
\delta_1^{\ast_s 2n} & = & {{2n}\choose{n}} 2^{-2n} \delta_0 + 2\sum_{k=1}^n {{2n}\choose{n-k}} 2^{-2n} \delta_{2k}, \\
\delta_1^{\ast_s (2n+ 1)} & = &2 \sum_{k=0}^n {{2n+1}\choose{n-k}} 2^{-2n-1} \delta_{2k+1}.
\end{eqnarray*}
Consequently we have
\begin{eqnarray*}
{\rm Exp}_{\ast_s} ( a \delta_1) & = & e^{-a} \sum_{n=0}^{\infty} \frac{1}{n!\, n!} \left(\frac{a}{2}\right)^{2n} \delta_0 + 2 e^{-a} \sum_{k=1}^{\infty} \sum_{n=0}^{\infty} \frac{\left({a/2}\right)^{2n+k}}{n!\, (n+k)!} \delta_{k} \\
& = & e^{-a} I_0(a) \delta_0 + 2 e^{-a} \sum_{k=1}^{\infty} I_k(a) \delta_{k},
\end{eqnarray*}
where $I_k$ is the modified Bessel function of the first kind.
\item[{\bf2a.}] In the same way we get that
\begin{eqnarray*}
\delta_1^{\ast_{s,\alpha} 2n} & = & {{2n}\choose{n}} 2^{-2n} \delta_0 + 2\sum_{k=1}^n {{2n}\choose{n-k}} 2^{-2n} \delta_{(2k)^{1/{\alpha}}}, \\
\delta_1^{\ast_{s,\alpha} (2n+ 1)} & = & 2 \sum_{k=0}^n {{2n+1}\choose{n-k}} 2^{-2n-1} \delta_{(2k+1)^{1/{\alpha}}},
\end{eqnarray*}
thus
\begin{eqnarray*}
{\rm Exp}_{\ast_{s,\alpha}} ( a \delta_1) & = & e^{-a} I_0(a) \delta_0 + 2 e^{-a} \sum_{k=1}^{\infty} I_k(a) \delta_{k^{1/{\alpha}}}.
\end{eqnarray*}
\item[{\bf 3a.}] For $\diamond = \ast_{\alpha}$, $\alpha>0$ we have
$$
\left( p \delta_1 + (1-p) \delta_0 \right)^{\star_{\alpha} n} = \sum_{k=0}^n {{n}\choose{k}} p^k (1-p )^{n-k} \delta_{k^{1/\alpha} }
$$
and
$$
{\rm Exp}_{\star_{\alpha}} ( a \delta_1) = e^{-a} \sum_{k=0}^{\infty} \frac{a^k}{k!} \;\delta_{k^{1/{\alpha}}}.
$$
Notice that $\star_{\alpha}$-generalized Poisson measure is purely atomic with the support equal $\{ 1, 2^{1/{\alpha}}, 3^{1/{\alpha}}, 4^{1/{\alpha}}, \dots \}$. Also in examples 1a and 2a the measures ${\rm Exp}_{\diamond} ( a \delta_1)$ are purely atomic. We see below that this property is exceptional for $\diamond$-generalized Poisson measure.
\item[{\bf 4a.}] For the Kendall $\vartriangle_{\alpha}$ generalized convolution, ${\alpha}>0$, we have
\begin{eqnarray*}
\left( p \delta_1 + (1-p) \delta_0 \right)^{\vartriangle_{\alpha} n}(dt) & = & (1-p)^n \, \delta_0(dt) \, + \,  np(1-p)^{n-1} \, \delta_1(dt) \\[1mm]
&&  \hspace{-35mm}+ \alpha p^2 n (n-1) t^{-2\alpha - 1} \left(1 - \frac{p}{t^{\alpha}}\right)^{n-2} \mathbf{1}_{(1,\infty)}(t)dt.
\end{eqnarray*}
To see this we denote by $F_n$ the cumulative distribution function of $\left( p \delta_1 + (1-p) \delta_0 \right)^{\vartriangle_{ \alpha} n}$ and calculate the corresponding generalized characteristic functions:
\begin{eqnarray*}
\lefteqn{ \int_0^{\infty} \bigl( 1 - (us)^{\alpha}\bigr)_{+} \left( p \delta_1 + (1-p) \delta_0\right)^{\vartriangle_{\alpha} n}(ds) } \\
 & = & \left[ \int_0^{\infty} \bigl( 1 - (us)^{\alpha}\bigr)_{+} \left( p \delta_1 + (1-p) \delta_0\right)(ds) \right]^n \\
 & = & \left[p \bigl( 1 - u^{\alpha}\bigr)_{+} + 1 -p\right]^n.
\end{eqnarray*}
On the other hand we have
\begin{eqnarray*}
\int_0^{\infty}\!\! \bigl( 1 - (us)^{\alpha}\bigr)_{+} dF_n(s)\!\! & = & \!\! \bigl( 1 - (us)^{\alpha}\bigr)_{+} F_n(s)\big|_0^{u^{-1}} + \frac{\alpha}{ u^{\alpha}} \int_0^{u^{-1}} \!\!\!\! s^{\alpha-1} F_n(s) ds  \\
& = &  - F(0) + \alpha u^{-\alpha} \int_0^{u^{-1}} \!\!\!\! s^{\alpha-1} F_n(s) ds.
\end{eqnarray*}
Substituting now $t=u^{-1}$ and comparing both formulations we obtain
$$
\int_0^t s^{\alpha - 1} F_n(s) ds = \left\{ \begin{array}{ll}
\alpha t^{-\alpha} \bigl(F(0) + 1 - p\bigr)^n & 0<t<1; \\
\alpha t^{-\alpha} \bigl(F(0) + (1 - p\, t^{-\alpha})^n \bigr)^n & t>1.
\end{array}\right.
$$
Differentiating this equality with respect to $t$ we have
$$
F_n(t)  = \left\{ \begin{array}{ll}
 \bigl(F(0) + 1 - p\bigr)^n & 0<t<1; \\
F(0) + \bigl(1 - p\, t^{-\alpha})^n \bigr)^{n-1}\bigl( 1 + (n-1)p\,t^{-\alpha}\bigr)  & t>1.
\end{array}\right.
$$
Since $\lim_{t\rightarrow \infty} F_n(t) = F(0) + 1$ we conclude $F(0)=0$. Now it is enough to notice that the function $F_n$ has two jumps at 0 and at 1 and it has an absolutely continuous part on $(1, \infty)$.
In a similar way we can show that
$$
 {\rm Exp}_{\vartriangle_{\alpha}}(a \delta_1) (du) = e^{-a} \delta_0 (du) + a e^{-a} \delta_1(du) + \frac{a^2 \alpha}{u^{2\alpha +1}} e^{-au^{-\alpha}} \mathbf{1}_{(1,\infty)}(u) du.
$$
We see that, except for the two atoms at zero and one, the $\vartriangle_{\alpha}$-generalized Poisson measure  is absolutely continuous with respect to the Lebesgue measure.
\item[{\bf 5a.}] For the technical reasons we consider here the special case of the generalized Kingman convolution  $\otimes_{\omega_3} \colon \mathcal{P}^2_{+} \rightarrow \mathcal{P}_{+}$. Since the generalized convolutions defined by $\omega_3$- the uniform distribution on the unit sphere in $\mathbb{R}^3$-  and convolution defined by its one-dimensional projection $\omega_{3,1}$ are the same and $\omega_{3,1}(du) = \frac{1}{2} \mathbf{1}_{[-1,1]}(u) du$ the calculations are simpler than in general case.
Since $\omega_{3,1}$ is the uniform distribution on $[-1,1]$ then $\omega_{3,1}^{\ast n}$ are also well known and e.g. in \cite{Kill} we can find that for $n \geqslant 1$ the measure $\omega_{3,1}^{\ast n}$, the classical $n$'th convolution of $\omega_{3,1}$, has  the following density function:
$$
f^{(n)}(x)\! =\! \left\{ \begin{array}{cl}
  {\displaystyle\!\!\!\!\sum_{i=0}^k}\!\left(-1\right)^i\!
  \!{{n}\choose{i}}\!\frac{\left(x+n-2i\right)^{n-1}}{\left(n-1\right)!\; 2^n}, & x \in [ -n + \!2k,-n\!+\!2\!\left(k+1\right)) \\[-1mm]
  & k=0, \ldots, n-1\\[2mm]
  0  &\!\! \hbox{{\rm otherwise}} \\
  \end{array}\!\!.\right.
$$
By the construction of Kingman generalized convolution we know that $\omega_{3,1}^{\ast n} = \delta_1^{\otimes_{\omega_{3,1} n}} \circ \omega_{3,1}$, thus the cumulative distribution function $F_n$ of the measure $\delta_1^{\otimes_{\omega_{3,1} n}}$ we obtain by solving the following equation:
$$
f_n(u) = \frac{1}{2} \int_{\mathbb{R}} y^{-1} \mathbf{1}_{[-1,1]}(u/y) dF_n(y).
$$
Since $F_n(0)=0$ then we can consider only $u>0$ and integrating by parts implies that
$$
f_n(u) =  - \frac{1}{2} u^{-1} F_n(u)  + \frac{1}{2} \int_u^{\infty} y^{-2} F_n(y) dy.
$$
Now it is enough to differentiate both sides of this equation to get
$$
F_n'(u) = - 2u f_n'(u)\mathbf{1}_{(0,\infty)}(u).
$$
Since $\int_0^{\infty} F_n'(u) du = 1$ we conclude that the measure $\delta_1^{\otimes_{\omega_{3,1} n}}$ is absolutely continuous with density $F_n'$.
These formulas can be used in calculating the measure $\Exp_{\otimes_{\omega_3}}(c \delta_1)$. However in
 \cite{reprLevy}, using much simpler method, it was shown that for any $c>0$
$$
{\rm Exp}_{\omega_{3,1}}(a \delta_1) = Exp \left(a\, \omega_{3,1} \right) \ast \left( \delta_0- a \,\omega_{3,1} + a\, \delta_1 \right),
$$
Notice that, except for two atoms, this measure contains an absolutely continuous with respect to the Lebesgue measure part.
\item[{\bf 6a.}] For the max-convolution we have
$$
\left( p \delta_0 + q \delta_1\right)^{{\bigcirc\hspace{-2.54mm}\vee \,\,} n} = p^n \delta_0 +(1-p^n) \delta_1; \quad {\rm Exp}_{{\bigcirc\hspace{-2.54mm}\vee \,\,}}(a \delta_1) = e^{-a} \delta_0 +(1-e^{-a}) \delta_1.
$$
\end{namelist}

\begin{defn}
A stochastic process $\{ N_I(t) \colon t \geqslant 0 \}$ is a $\diamond$-generalized Poisson process of type I  with intensity $c>0$, where $\diamond$ is a generalized convolution, if it is a Markov process with the transition probability
$$
P_{s,t}(x,\cdot\,)=\delta_{x}\diamond{\rm Exp}_{\diamond}(c (t-s) \delta_1)(\cdot\,),\quad x\in \mathbb{R}_{+},\: s<t
$$
and $\mathbf{P}\{ N_I(0) =0\} = 1$.
\end{defn}

The consistency of this definition, proof that $P_{s,t}(x,\cdot\,)$ form a consistent family of transition probabilities for which the Chapman-Kolmogorov equations hold one can find in \cite{WeakLevyProc}.

Of course we have that $N_I(t)$ has the distribution ${\rm Exp}_{\diamond}(c t \delta_1)$. Moreover this process has independent increments in the sense of generalized convolution, i.e. for every $t>s$ there exists a random variable $X_{[s,t)}$ with distribution ${\rm Exp}_{\diamond}(c (t-s) \delta_1)$ independent of $N_I(s)$ such that
$$
\mathcal{L}(N_I(t)) = \mathcal{L}(N_I(s)) \diamond \mathcal{L}(X_{[s,t)}),
$$
since
$$
{\rm Exp}_{\diamond}(c t \delta_1) = {\rm Exp}_{\diamond}(c (t-s) \delta_1) \diamond {\rm Exp}_{\diamond}(c s \delta_1).
$$
The variables $X_{[s,t)}$ for $0<s<t$ are playing the same role as increments in the classical Poisson process, however they are determined here only in the sense of uniqueness of their distributions. Similarly like in the classical Poisson process their distributions are stationary in time.

The existence of the $\diamond$-generalized Poisson process of type I follows from the existence of the Markov process with given transition probabilities and from infinite divisibility of ${\rm Exp}_{\diamond}(c t \delta_1)$ in the sense of generalized convolution.

This kind of $\diamond$-generalized Poisson process usually cannot be used for counting buses coming to the station since the distribution of $N_I(t)$ is not taking values in the set of natural numbers. In the case of Kendall convolution we could eventually model primitive counting: zero-one-many, which is not very interesting.

The next proposition describes the infinitesimal operators $A$ for $\diamond$-generalized Poisson processes type I, where
$$
A f(x) \stackrel{def}{=} \lim_{s\rightarrow t^-} \frac{\mathbf{E}\left(f(N_I(t)) \big| N_I(s) = x\right) - f(x)}{t-s}.
$$
\begin{prop}
Let $f$ be a continuous bounded function on $\mathbb{K}$ and let $\{N_I(t) \colon t \geqslant 0 \}$ be a $\diamond$-generalized Poisson process type I with intensity $c>0$, where $\diamond$ is a generalized convolution. Then
$$
Af(x) = c \int_{\mathbb{K}} \left( f(r) - f(x) \right) \delta_x \diamond \delta_1(dr).
$$
\end{prop}

\noindent
{\bf Proof.} Since the distribution of $N_I(t)$ given $N_I(s) = x$ is the transition probability $\delta_x \diamond \Exp(c(t-s) \delta_1)$ we see that
\begin{eqnarray*}
\lefteqn{\hspace{-17mm}\mathbf{E}\left(f(N_I(t)) \big| N_I(s) = x\right) = \int_{\mathbb{K}}f(r) \delta_x \diamond \Exp\left(c(t-s)\delta_1\right)(dr)} \\
& = & e^{-c(t-s)} \sum_{k=0}^{\infty} \frac{(c(t-s))^k}{k!} \int_{\mathbb{K}}f(r) \delta_x \diamond \delta_1^{\diamond k} (dr) \\
& = & e^{-c(t-s)} f(x) +  c(t-s) \int_{\mathbb{K}}f(r) \delta_x \diamond \delta_1 (dr) \\
& + & e^{-c(t-s)} \sum_{k=2}^{\infty} \frac{(c(t-s))^k}{k!} \int_{\mathbb{K}}f(r) \delta_x \diamond \delta_1^{\diamond k} (dr). \end{eqnarray*}
Since $f$ is bounded, let say $|f(x)| < M$ for every $x \in \mathbb{K}$, and $ \delta_x \diamond \delta_1^{\diamond k}$ is a probability measure then for $s\rightarrow t^-$
$$
\left| \sum_{k=2}^{\infty} \frac{(c(t-s))^{k-1}}{k!} \int_{\mathbb{K}}f(r) \delta_x \diamond \delta_1^{\diamond k} (dr)\right| \leqslant M \sum_{k=2}^{\infty} \frac{(c(t-s))^{k-1}}{k!} \rightarrow 0.
$$
Consequently
\begin{eqnarray*}
Af(x) & = & \lim_{s\rightarrow t^-} \left[ \frac{e^{-c(t-s)} - 1}{t-s} f(x) + c \int_{\mathbb{K}}f(r) \delta_x \diamond \delta_1 (dr) e^{-c(t-s)} \right] \\
 & = & c \int_{\mathbb{K}} \left( f(r) - f(x) \right) \delta_x \diamond \delta_1(dr).
\end{eqnarray*}
\qed

\vspace{2mm}

\noindent {\bf Examples.}
\begin{namelist}{ll}
\item[{\bf 1b.}] If $\{ N_I(t) \colon t \geqslant 0 \}$ is the $\ast_s$-generalized Poisson process of the first kind with respect to the symmetric convolution $\ast_s$ then the transition probability is given by
    \begin{eqnarray*}
    \delta_x \ast_s {\rm Exp}_{\ast_s} (c(t-s)\delta_1) & = &   e^{-c(t-s)} I_0(c(t-s)) \delta_x \\
    & + & e^{-c(t-s)} \sum_{k=1}^{\infty} I_k(c(t-s)) \left( \delta_{|x-k|} + \delta_{x+k} \right)
    \end{eqnarray*}
    and the infinitesimal operator $A$ for this process on the class of bounded measurable functions takes the form
    $$
    A f(x) =c \int_0^{\infty} \bigl( f(r) - f(x) \bigr) \delta_x \ast_s \delta_1 (dr) = c \bigl( f(1+x) + f(|1-x|) - f(x)\bigr).
    $$
\item[{\bf 2b.}] If $\{ N_I(t) \colon t \geqslant 0 \}$ be the $\ast_{s,\alpha}$-generalized Poisson process of the first kind with respect to the $(\alpha,1)$- convolution $\ast_{s,\alpha}$ then the transition probability is given by
    \begin{eqnarray*}
    \delta_x \ast_s {\rm Exp}_{\ast_{s,\alpha}} (c(t-s)\delta_1) & = &   e^{-c(t-s)} I_0(c(t-s)) \delta_x \\
    & + & e^{-c(t-s)} \sum_{k=1}^{\infty} I_k(c(t-s)) \left( \delta_{|x-k|^{\alpha}} + \delta_{(x+k)^{\alpha}} \right)
    \end{eqnarray*}
    and the infinitesimal operator $A$ for this process on the class of bounded measurable functions takes the form
    $$
    A f(x) = c\!\int_0^{\infty}\!\!\!\! \bigl( f(r) - f(x) \bigr) \delta_x \ast_{s,\alpha} \delta_1 (dr) = c \bigl( f((1+x)^{\alpha}) + f(|1-x|^{\alpha}) - f(x)\bigr).
    $$
\item[{\bf 3b.}] If $\{ N_I(t) \colon t \geqslant 0 \}$ is the $\ast_{\alpha}$-generalized Poisson process of type I, $\alpha >0$, then
    $$
    \delta_x \ast_{\alpha} {\rm Exp}_{\ast_{\alpha}} (c(t-s)\delta_1) = e^{-c(t-s)} \sum_{k=0}^{\infty} \frac{\left(c(t-s)\right)^k}{k!} \delta_{(x^{\alpha} +k)^{1/{\alpha}}},
    $$
    and
    $$
    Af(x) = c \int_0^{\infty} \left( f(r) - f(x) \right) \delta_x \ast_{\alpha} \delta_1 (dr) = c \left( f\left( (x^{\alpha} +1)^{1/{\alpha}} \right) - f(x) \right),
    $$
    for every bounded measurable function $f$ on $[0,\infty)$.
\item[{\bf 4b.}] Let $\{ N_I(t) \colon t \geqslant 0 \}$ be the $\vartriangle_{\alpha}$-generalized Poisson process of type I with respect to the Kendall convolution, $\alpha >0$ with intensity $c>0$. The easiest way to calculate the distribution of the transition probability  $\nu_x = \delta_x \vartriangle_{\alpha}{\rm Exp}_{\vartriangle_{\alpha}}(c (t-s) \delta_1)$ with the cumulative distribution function $F_x$ lies in calculating the corresponding generalized characteristic functions:
$$
\int_0^{t^{-1}} \bigl( 1 - (ts)^{\alpha}\bigr)_{+} dF_x(s) = \bigl( 1 - (xt)^{\alpha}\bigr)_{+} \exp \bigl\{ - c(t-s) \bigl( t^{\alpha} \wedge 1\bigr)\bigr\}.
$$
Integrating by parts the integral on the left hand side and then substituting $t = u^{-1}$ we obtain
$$
\int_0^u s^{\alpha-1} F(s) ds = \alpha^{-1} u^{\alpha} \bigl( 1 - u^{-\alpha} s^{\alpha}\bigr)_{+} \exp \bigl\{ - c(t-s) \bigl( u^{-\alpha} \wedge 1\bigr)\bigr\}.
$$
Now it is enough to differentiate both sides of this equality with respect to $u$ to obtain
$$
F_x(u) = \left\{ \begin{array}{lll}
   \left( 1 + \frac{c z}{u^{\alpha}} - \frac{c z x^{\alpha}} {u^{2\alpha}} \right) e^{- c z u^{-\alpha}}\mathbf{1}_{(x,\infty)} & \hbox{if} &\!\! x \geqslant 1; \\[2mm]
   e^{-cz} \mathbf{1}_{(x,1]}(u) + \left( 1 + \frac{c z}{u^{\alpha}} - \frac{c z x^{\alpha}}{u^{2\alpha}} \right) e^{- c z u^{-\alpha}}\mathbf{1}_{(1,\infty)}(u) & \hbox{if} &\!\! x < 1,
   \end{array} \right.
$$
where $z=t-s$.
By Proposition 2, the infinitesimal operator $A$ of the process $\{N_I(t) \colon t \geqslant 0 \}$ for real measurable bonded function $f$ on $[0,\infty)$ is given by:
$$
Af(x) =  \left\{\!\! \begin{array}{lll}
  - c x^{-\alpha} f(x) + 2 \alpha c x^{\alpha} \int_x^{\infty} f(r) r^{-2\alpha -1} dr & \hbox{if} & x \geqslant 1; \\[2mm]
  - c f(x) + c\bigl(1-x^{\alpha}\bigr) + 2 \alpha c x^{\alpha} \int_1^{\infty} f(r) r^{-2\alpha -1} dr & \hbox{if} & x < 1.
  \end{array}\right.
$$
\item[{\bf 6b.}] If $\{ N_I(t) \colon t \geqslant 0 \}$ is the $\bigcirc\hspace{-3.54mm}\vee \,\,$-generalized Poisson process of the first kind with respect to the max-convolution then
    $$
    \delta_x \Max {\rm Exp}_{{\bigcirc\hspace{-2.54mm}\vee \,\,}}(c(t-s) \delta_1) = \left\{ \begin{array}{lll}
    e^{-c(t-s)} \delta_1 + (1- e^{-c(t-s)} ) \delta_x & \hbox{if} & x\in[0,1] \\
    \delta_x  & \hbox{if} & x > 1,
    \end{array}\right.
    $$
    and the infinitesimal operator
    $$
    Af(x) = \left\{ \begin{array}{ll}
     \!c \left( f(1) - f(x) \right) & x \in [0,1], \\
     \!0 & x>1. \end{array}\right.
    $$
\end{namelist}

\section{Monotonicity property  and random walk under generalized convolution}
For the classical Poisson process the corresponding process of change times $T_1 + \dots T_n$, $n\in \mathbb{N}$, is also monotonically increasing. In order to get this property we will concentrate on monotonic generalized convolutions on $\mathcal{P}_{+}$, where
\begin{defn}
A generalized convolution $\diamond$ is monotonic if for all $x,y \geqslant 0$, $x\vee y =\max\{x,y\}$
$$
\delta_x \diamond \delta_y \left( [x\vee y, \infty) \right) = 1.
$$
\end{defn}
It is easy to see that classical convolution, stable $\ast_{\alpha}$ convolution and the Kendall $\vartriangle_{\alpha}$ convolution are monotonic. The best known Kingman convolution $\otimes_{\omega_s}$ is not monotonic since the measure $\delta_x \otimes_{\omega_s} \delta_y$ has the support equal $[|x-y|, x+y]$.
\begin{lem}
Assume that $\diamond$ is a monotonic generalized convolution. If $X$ is a nonnegative random variable and $x\diamond y$, $x,y \geqslant 0$, is any random variable with distribution $\delta_x \diamond \delta_y$ independent of $X$ then
$$
\mathbf{P} \left\{ X > x \diamond y, X>x \right\} = \mathbf{P} \left\{ X > x \diamond y \right\}.
$$
\end{lem}

{\bf Proof.} By monotonicity definition the measure $\delta_x \diamond \delta_y$ lives on the set $[x, \infty)$. Thus we have
\begin{eqnarray*}
\lefteqn{\mathbf{P} \left\{ X > x \diamond y, X>x \right\} = \int_0^{\infty}  \mathbf{P} \left\{ X > u , X>x \right\} \delta_x \diamond \delta_y (du)} \\
 & = & \int_x^{\infty}  \mathbf{P} \left\{ X > u , X>x \right\} \delta_x \diamond \delta_y (du) \\
 & = & \int_x^{\infty}  \mathbf{P} \left\{ X > u  \right\} \delta_x \diamond \delta_y (du) = \mathbf{P} \left\{ X > x \diamond y \right\}.
\end{eqnarray*}
\qed

 \begin{defn}
Let $\diamond$ be a generalized convolution on $\mathcal{P}_{+}$ and $\{T_n \colon n \in \mathbb{N} \}$ be a sequence of independent identically distributed random variables with the distribution $\nu \in \mathcal{P}_{+}$. The Markov process $(S_n)$ with $S_1 = T_1$, the transition probabilities
$$
P_{k,n}(x,\cdot )=P(S_n\in \cdot\,|S_k=x):= \delta_x \diamond \nu^{\diamond (n-k)}(\cdot\,)
$$
such that $ S_n$  is independent of the variables $T_{n+1}, T_{n+2}, \dots$ we call a generalized random walk with respect to the convolution $\diamond$ with the step distribution $\nu$.
\end{defn}
The Kendall random walk was studied by Jasiulis-Go{\l}dyn in \cite{KendallWalk}.
We present here few basic properties of the random walk under monotonic generalized convolution; these properties will be used in the next section for the Poisson process of  type II. By $F_n$ we denote the cumulative distribution function of the variable $S_n$.

\begin{lem}\label{rem:0}
If the generalized convolution $\diamond$ is monotonic then for every $k_1< k_2< \dots< k_{n+1}$
\begin{eqnarray*}
{\rm a)}  & \mathbf{P}\left\{ S_{k_{n+1}} > x, S_{k_n} > x , \dots, S_{k_1}> x\right\} = 1-F_{k_1}(x), \\
{\rm b)} & \mathbf{P}\left\{ S_{k_{n+1}} > x, S_{k_n} \leqslant x, \dots, S_{k_1} \leqslant x \right\} = \mathbf{P}\left\{ S_{k_{n+1}} > x, S_{k_n} \leqslant x \right\} \\
 & = F_{k_n}(x) - F_{k_{n+1}}(x).
\end{eqnarray*}
\end{lem}

{\bf Proof.} Notice first that for every $w>x$ we have
$$
\delta_w \diamond \delta_y \bigl((x, \infty) \bigr) \geqslant \delta_w \diamond \delta_y \bigl((w, \infty) \bigr) = 1.
$$
In order to explain method of calculations it is enough to concentrate on three variables. For $n,k,\ell \in \mathbb{N}$ we have
\begin{eqnarray*}
\lefteqn{\mathbf{P} \left\{ S_{n+k+\ell} > x, S_{n+k}> x, S_n > x \right\} } \\
 & = & \int_x^{\infty} \int_0^{\infty} \int_0^{\infty} \delta_u \diamond \delta_v \diamond \delta_w \bigl((x, \infty) \bigr) dF_{\ell}(u) dF_k(v) dF_n(w) \\
  & = & \int_x^{\infty} \int_0^{\infty} \int_0^{\infty} \int_0^{\infty}\delta_z \diamond \delta_w \bigl((x, \infty) \bigr) \delta_u \diamond \delta_v (dz) dF_{\ell}(u) dF_k(v) dF_n(w) \\
 & = & \int_x^{\infty} \int_0^{\infty} \int_0^{\infty} \int_0^{\infty} \delta_u \diamond \delta_v (dz) dF_{\ell}(u) dF_k(v) dF_n(w)  \\
& = & \int_x^{\infty} \int_0^{\infty} \int_0^{\infty} dF_{\ell}(u) dF_k(v) dF_n(w) = 1 - F_n(x).
\end{eqnarray*}
The third equality holds because $w>x$, the forth one, because $\delta_u \diamond \delta_v$ is a probability measure. Consequently we have
\begin{eqnarray*}
\lefteqn{\mathbf{P} \left\{ S_{n+k+\ell} > x, S_{n+k}\leqslant x, S_n\leqslant x \right\} } \\
 & = & \mathbf{P} \left\{ S_{n+k+\ell} > x, S_{n+k}\leqslant x \right\} - \mathbf{P} \left\{ S_{n+k+\ell} > x, S_{n+k}\leqslant x, S_n> x \right\} \\
 & = & \mathbf{P} \left\{ S_{n+k+\ell} > x \right\} - \mathbf{P} \left\{ S_{n+k+\ell} > x, S_{n+k}> x \right\} \\
 & - & \mathbf{P} \left\{ S_{n+k+\ell} > x, S_{n}\leqslant x \right\} - \mathbf{P} \left\{ S_{n+k+\ell} > x, S_{n+k}> x, S_n> x \right\} \\
& = & 1 - F_{n+k+\ell}(x) - (1-F_{n+k}(x)) -(1-F_n(x)) + 1-F_n(x) \\
 & = & F_{n+k}(x) - F_{n+k+\ell}(x).
\end{eqnarray*}
\qed

\begin{lem}\label{lem:6}
Let $(\mathcal{P},\diamond)$ be monotonic convolution algebra. For each choice of\,  $0\leqslant s \leqslant t$ and for every $k,n\in \mathbb{N}_0$ we have
\begin{eqnarray*}
\lefteqn{\hspace{-10mm}\mathbf{P} \left\{ S_{n+k+1} > t, S_{n+k} \leqslant t , \dots,  S_{k+2} \leqslant t , S_{k+1} > s, S_k \leqslant s,  \cdots , S_1 \leqslant s, \right\}} \\
& = & \mathbf{P} \left\{ S_{n+k+1} > t, S_{n+k} \leqslant t , S_{k+1} > s, S_k \leqslant s \right\} \\
& = & \mathbf{P} \left\{ S_{n+k}\leqslant t, S_k \leqslant s \right\} - \mathbf{P} \left\{ S_{n+k} \leqslant t, S_{k+1} \leqslant s \right\} \\
& - & \mathbf{P} \left\{ S_{n+k+1} \leqslant t, S_k \leqslant s  \right\} + \mathbf{P} \left\{ S_{n+k+1} \leqslant t, S_{k+1} \leqslant s \right\}.
\end{eqnarray*}
\end{lem}

{\bf Proof.}
Monotonicity property implies that:
\begin{eqnarray*}
\lefteqn{\mathbf{P} \left\{ S_{n+k+1} > t, S_{n+k} \leqslant t , \dots,  S_{k+2} \leqslant t , S_{k+1} > s, S_k \leqslant s,  \dots , S_1 \leqslant s, \right\}} \\
& = & \mathbf{P} \left\{ S_{n+k} \leqslant t, \dots,  S_{k+2} \leqslant t , S_{k+1} > s, S_k \leqslant s,  \dots, S_1 \leqslant s, \right\} \\
& - & \mathbf{P} \left\{ S_{n+k+1} \leqslant t,  \dots,  S_{k+2}< t , S_{k+1} > s, S_k \leqslant s,  \dots , S_1 \leqslant s, \right\} \\
& = & \mathbf{P} \left\{ S_{n+k} \leqslant t, \dots,  S_{k+2} \leqslant t , S_{k+1} \leqslant t, S_k \leqslant s,  \dots , S_1 \leqslant s, \right\} \\
& - & \mathbf{P} \left\{ S_{n+k} \leqslant t, \dots,  S_{k+2} \leqslant t , S_{k+1} \leqslant s, S_k \leqslant s,  \dots , S_1 \leqslant s, \right\} \\
& - & \mathbf{P} \left\{ S_{n+k+1} \leqslant t, \dots,  S_{k+2}\leqslant t , S_{k+1} \leqslant t, S_k \leqslant s,  \dots , S_1 \leqslant s, \right\} \\
& + & \mathbf{P} \left\{ S_{n+k+1} \leqslant t, \dots,  S_{k+2}\leqslant t , S_{k+1} \leqslant s, S_k \leqslant s,  \dots , S_1 \leqslant s, \right\} \\
& = & \mathbf{P} \left\{ S_{n+k}\leqslant t, S_k \leqslant s \right\} - \mathbf{P} \left\{ S_{n+k} \leqslant t, S_{k+1} \leqslant s \right\} \\
& - & \mathbf{P} \left\{ S_{n+k+1} \leqslant t, S_k \leqslant s  \right\} + \mathbf{P} \left\{ S_{n+k+1} \leqslant t, S_{k+1} \leqslant s \right\}
\end{eqnarray*}
\qed
\section{Lack of the memory property}
The classical direct construction of the Poisson process $\{ N(t) \colon t \geqslant 0\}$ is based on the sequence $(T_k)$ of i.i.d. random variables with the exponential distribution $\Gamma(1,a)$.
In fact the main reason why the exponential distribution is used in this construction is lack of the memory, the exclusive property for exponential distribution. It turns out that this property strongly depends on the considered convolution.

\begin{defn}
A probability distribution $\nu$ with the cumulative distribution function $F$ has the lack of the memory property with respect to generalized convolution $\diamond$ if
$$
\mathbf{P}\bigl\{ X > x \diamond y  \big| X > x \bigr\} =  1 - F(y), \quad x,y \in \mathbb{R}_{+},
$$
where $X$ with distribution function $F$ is independent of the random variable $x \diamond y$ having distribution $\delta_x \diamond \delta_y$.
\end{defn}

\vspace{2mm}

\begin{prop}\label{monoton}
The distribution function $F\not\equiv \mathbf{1}_{[0,\infty)}$ has lack of the memory property with respect to the monotonic generalized convolution $\diamond$ if and only if the algebra $(\mathcal{P}_{+}, \diamond)$ is regular with homomorphism $h(\delta_t)$, $t \geqslant 0$, which is monotonically decreasing as a function of $t$ and $F(t) = 1 - h(\delta_{c^{-1}t})$ for some $c>0$.
\end{prop}

\noindent
{\bf Proof.}
For the monotonic generalized convolution by Lemma 1 the lack of memory property condition takes the form
$$
\mathbf{P}\bigl\{ X > x \diamond y \bigr\} = \bigl(1 - F(x)\bigr) \bigl(1 - F(y)\bigr), \quad x,y \in \mathbb{R}_{+}.
$$
Assume first that $h(\delta_t) = 1- F(ct)$ is a nontrivial homomorphism for $\diamond$. Since $h(\delta_0 = 1$, then $F(0) =0$. Since $h$ is continuous and monotonically decreasing then $F$ is continuous and increasing on $[0,\infty)$. Since the convolution $\diamond$ is monotonic then $\delta_x \diamond \delta_x$ as a measure concentrated on $[x,\infty)$ converges weakly to $\delta_{\infty}$ when $x \rightarrow \infty$. Consequently $h(\delta_{\infty}) h(\delta_{\infty}) = h(\delta_{\infty}) = 0$ and $\lim_{x\rightarrow \infty} F(x) = 1$. Moreover
$$
h(\delta_x \diamond \delta_y ) = \left( 1 - F(x) \right) \left( 1 - F(y) \right).
$$
On the other hand
\begin{eqnarray*}
h(\delta_x \diamond \delta_y ) & = & \int_0^{\infty} h(\delta_z) (\delta_x \diamond \delta_y) (dz) = \int_0^{\infty} (1 - F(z)) (\delta_x \diamond \delta_y) (dz) \\
 & = & \int_0^{\infty} \mathbf{P} \{ X > z \}( \delta_x \diamond \delta_y) (dz) = \mathbf{P} \{ X > x \diamond y\},
\end{eqnarray*}
thus $F$ is a cumulative distribution function with lack of memory property.

Assume now that $F$ has lack of memory property and define
$$
h(\lambda):= \mathbf{P} \{ X > \theta\} = \int_0^{\infty} (1-F(x)) \lambda(dx) = \int_0^{\infty} h(\delta_x)  \lambda(dx),
$$
where $\theta$ is a random variable with distribution $\lambda$ independent of $X$. By assumption $h(\delta_x \diamond \delta_y ) = h(\delta_x ) h(\delta_y)$. Since
\begin{eqnarray*}
h(\lambda_1 \diamond \lambda_2) & = & \int_0^{\infty} h(\delta_z) \lambda_1 \diamond \lambda_2 (dz) \\
& = & \int_0^{\infty} \int_0^{\infty} h(\delta_x \diamond \delta_y) \lambda_1(dx) \lambda_2(dy) \\
& = & \int_0^{\infty} \int_0^{\infty} h(\delta_x )h(\delta_y) \lambda_1(dx) \lambda_2(dy) = h(\lambda_1) h(\lambda_2)
\end{eqnarray*}
and for all $\lambda_1, \lambda_2 \in \mathcal{P}_{+}$, $p \in [0,1]$, $q=1-p$
$$
h(p \lambda_1 + q \lambda_2) = \int_0^{\infty} h(\delta_z) (p \lambda_1 + q \lambda_2)(dz) = p h(\lambda_1) + q h(\lambda_2).
$$
Finally, since $F \not\equiv 1$ we obtain that $h$ is a homomorphism for $\diamond$. \qed

\vspace{2mm}

\noindent 
\begin{rem}
{\rm  In particular it follows from Proposition \ref{monoton} that the regular generalized convolution $\diamond$ admits the existence of cumulative distribution function $F$ with the lack of memory property with respect to $\diamond$ if and only if the corresponding homomorphism $h$ is the tail of some distribution function. The uniqueness (up to a scale) of the homomorphism $h$ on the algebra $(\mathcal{P}_{+}, \diamond)$ implies also that each regular generalized convolution admits the existence of at most one (up to a scale) distribution with the lack of the memory property. }
\end{rem}

\noindent 
\begin{rem}
{\rm It is evident now that stable convolution and Kendall convolution admit existence of the distribution with the lack of the memory property since their homomorphisms $h(\delta_t) = e^{-t^{\alpha}}$ and $h(\delta_t) = (1-t^{\alpha})_{+}$ respectively are the tails of some distribution functions. Much more complicated seems to be checking whether the given cumulative distribution function $F$ can be distribution with the lack of memory property with respect to some generalized convolution. We see that this is equivalent with checking if $h = 1 - F$ can be a probability kernel for some generalized convolution. For the cumulative distribution function $F_{\Gamma}$ given by
$$
1 - F_{\Gamma}(t) = \frac{\Gamma(a,t)}{\Gamma(a)} = \frac{1}{\Gamma(a)} \int_t^{\infty} x^{a-1} e^{-x} dx, \quad \quad t>0,
$$
it is known that $h = 1 - F_{\Gamma}$ is the homomorphism for Kucharczak convolution (for details see \cite{KU2}). However for the following the cumulative distribution function
$$
F_{\beta}(t) = \frac{\Gamma(a+b)}{\Gamma(a) \Gamma(b)} \int_0^t x^{a-1} (1-x)^{b-1}, \quad \quad x \in (0,1)
$$
the positive answer for this question is equivalent with the existence for all $r,s >0$ a probability measure $\lambda_{r,s}$ such that
$$
\bigl( 1 - F_{\beta}(rt) \bigr) \bigl( 1 - F_{\beta}(st) \bigr) = \int_0^{\infty} \bigl( 1 - F_{\beta}(xt) \bigr) \lambda_{r,s} (dx).
$$
We do not know whether such measure $\lambda_{r,s}$ exists or not.}
\end{rem}

\section{$\diamond$-Generalized Poisson process of type II}
\begin{defn}
Let $\diamond$ be a generalized convolution, let $(T_k)$ be a sequence of i.i.d. random variables with distribution $\nu$, and $(S_n)$ be the corresponding random walk. Then the process $\{ N_{II}(t) \colon t \geqslant 0\}$ given by
$$
N_{II}(t) = \left\{ \begin{array}{l}
 \inf\{ n \colon S_{n+1} > t \} \\[2mm]
 \infty \; \hbox{ if such $n$ does not exists }.
 \end{array} \right.
$$
is called the $\diamond$-generalized Poisson process of type II. The $\diamond$-generalized Poisson process of type II is proper if $\nu$ has the lack of memory property.
\end{defn}

\vspace{2mm}

\noindent
{\bf Property 1.} For every generalized convolution $\diamond$
$$
\mathbf{P}\{ N_{II}(t) = 0 \} = \nu( (t,\infty)= 1 - F(t).
$$

\vspace{2mm}

\noindent
{\bf Property 2.} For every monotonic generalized convolution $\diamond$
$$
\mathbf{P}\{ N_{II}(t) = n \} = F_{n}(t) - F_{n+1}(t),
$$
where $F_n$ is the cumulative distribution function of $S_n$. Consequently
$$
\mathbf{E}N_{II}(t) = \sum_{n=1}^{\infty} F_n(t), \quad \mathbf{E}N_{II}^2(t) = 2 \sum_{n=1}^{\infty} n F_n(t) - \sum_{n=1}^{\infty} F_n(t).
$$

\vspace{2mm}

\noindent
{\bf Property 3.}
For every monotonic generalized convolution by Lemma \ref{lem:6} we have for every $0 \leqslant s < t$
\begin{eqnarray*}
\lefteqn{\mathbf{P} \left\{ N_{II}(t) - N_{II}(s) = n, \, N_{II}(s) = k \right\}} \\
 &= & \mathbf{P} \left\{ S_{n+k+1} > t, S_{n+k} \leqslant t , S_{k+1} > s, S_k \leqslant s \right\} \\
& = & \mathbf{P} \left\{ S_{n+k}\leqslant t, S_k \leqslant s \right\} - \mathbf{P} \left\{ S_{n+k} \leqslant t, S_{k+1} \leqslant s \right\} \\
& - & \mathbf{P} \left\{ S_{n+k+1} \leqslant t, S_k \leqslant s  \right\} + \mathbf{P} \left\{ S_{n+k+1} \leqslant t, S_{k+1} \leqslant s \right\}.
\end{eqnarray*}

\subsection{Proper $\star_{\alpha}$- generalized Poisson process of type II}

The function $\overline{F}(t) = e^{-t^{\alpha}}$, $\alpha > 0$,  is a homomorphism on the algebra $(\mathcal{P}_{+}, \ast_{\alpha})$ and it can be treated as the tail of some distribution on $[0,\infty)$. Consequently the Weibull distribution with the cumulative distribution function $F(t) = 1 - e^{-t^{\alpha}}$, $t>0$, has the lack of memory property with respect to $\ast_{\alpha}$ generalized convolution. In particular for $\alpha = 1$ we obtain that the exponential distribution $F(t) = 1 - e^{-t}$ has the lack of memory property with respect to the classical convolution.

Let $\nu$ with the distribution function $F$ be the step distribution in our construction. The generalized characteristic function of $S_1$ is given by
$$
\Phi_{\nu} (t) := h\left( T_t \nu\right) = \int_0^{\infty}\!\! h(\delta_{st}) dF(s) =  \int_0^{\infty}\!\! e^{- t^{\alpha} s^{\alpha}} \alpha s^{\alpha -1} e^{- s^{\alpha}} ds = \frac{1}{ t^{\alpha} + 1}.
$$
Consequently, the generalized characteristic function of $S_n$ with the distribution $\nu^{\ast_{\alpha} n}$ is given by
$$
\Phi_{\nu^{\ast_{\alpha} n}} (t) = \left( \Phi_{\nu} (t)\right)^n = \left(  t^{\alpha} + 1 \right)^{-n} = \int_0^{\infty} e^{- t^{\alpha} s^{\alpha}} \frac{\alpha}{\Gamma(n)} \, s^{\alpha n -1} e^{- s^{\alpha}} ds,
$$
which means that the distribution function $F_n = F^{\ast_{\alpha} n}$ has the density
$$
f_n(s) = \frac{\alpha}{\Gamma(n)}\, s^{\alpha n -1} e^{- s^{\alpha}} \mathbf{1}_{(0,\infty)}(s).
$$
The stable convolution can be easily written in the language of random variables:
$$
\mathcal{L}(X) \ast_{\alpha} \mathcal{L} (Y) = \mathcal{L} \left(\left(X^{\alpha} + Y^{\alpha}\right)^{1/{\alpha}} \right),
$$
for all nonnegative independent random variables $X$ and $Y$. Consequently the sequence $\{ S_n \colon n \in \mathbb{N}_0 \}$ can be written as:
$$
S_0 \equiv 0, \quad S_n = \left( T_1^{\alpha} + \dots +T_n^{\alpha}\right)^{1/{\alpha}},
$$
where $T_n$, $n \in \mathbb{N}$ is a sequence of i.i.d. random variables with the cumulative distribution function $F$. Using this representation and Property 3 we can calculate for all $t > s >0$
\begin{eqnarray*}
\lefteqn{ \mathbf{P} \left\{ N_{II}(t) = n+k, \, N_{II}(s) = k \right\} =} \\
 && \mathbf{P} \left\{ S_{n+k+1} > t, S_{n+k} \leqslant t, S_{k+1} > s, S_k \leqslant s \right\} = \\
&&  \hspace{-5mm} \int_0^s f_k(x) \int_{(s^{\alpha} - x^{\alpha})_{+}^{1/{\alpha}}}^{\infty} \hspace{-3mm}f_1(y) \int_0^{(t^{\alpha} - x^{\alpha} - y^{\alpha} )_{+}^{1/{\alpha}}} \hspace{-5mm} f_{n-1}(z) \int_{(t^{\alpha} - x^{\alpha} - y^{\alpha} -z^{\alpha})_{+}^{1/{\alpha}}}^{\infty}  \hspace{-18mm} f_1(u) du\, dz\, dy\, dx \\
 && = \frac{(t^{\alpha} - s^{\alpha})^n}{n!}\, \frac{s^{\alpha k}}{k!}\,\, e^{- t^{\alpha}}.
\end{eqnarray*}
We obtain that, similarly to the classical Poisson process, the process  $N_{II}$ has independent increments, since
\begin{eqnarray*}
\lefteqn{\hspace{-3mm}\mathbf{P} \left\{ N_{II}(t)- N_{II}(s) = n \right\} = \sum_{k=0}^{\infty} \mathbf{P} \left\{ N_{II}(t) = n+k, \, N_{II}(s) = k \right\} }\\
& = & \sum_{k=0}^{\infty} \frac{(t^{\alpha} - s^{\alpha})^n}{n!}\, \frac{s^{\alpha k}}{k!}\,\, e^{- t^{\alpha}} = \frac{(t^{\alpha} - s^{\alpha})^n}{n!}\,\, e^{- ( t^{\alpha} - s^{\alpha})}
\end{eqnarray*}
and
\begin{eqnarray*}
\mathbf{P} \left\{ N_{II}(s) = k \right\} & = & \int_0^s \frac{\alpha}{\Gamma(k)} \, y^{\alpha k -1} e^{- y^{\alpha}} \int_{(s^{\alpha} - y^{\alpha})^{1/{\alpha}} } \alpha x^{\alpha -1} e^{-x^{\alpha}} \, dx\, dy \\
& = &  \frac{\alpha}{\Gamma(k)}\, e^{-s^{\alpha}} \int_0^s y^{\alpha k -1} \, dy =  \frac{s^{\alpha k}}{k!}\, e^{-s^{\alpha}}.
\end{eqnarray*}
This implies that
\begin{eqnarray*}
\lefteqn{\mathbf{P} \left\{ N_{II}(t)- N_{II}(s) = n,\, N_{II}(s) = k \right\} } \\
  & = & \mathbf{P} \left\{ N_{II}(t)- N_{II}(s) = n \right\} \mathbf{P} \left\{ N_{II}(s) = k \right\}.
\end{eqnarray*}
In a sense the increments of this process are stationary if we replace the Lebesgue measure on $[0,\infty)$ by the set function defined by $\ell([a, b)) = (b^{\alpha} - a^{\alpha})^{1/{\alpha}}$ since for $t>s$
$$
N_{II}(t)- N_{II}(s) \stackrel{d}{=} N_{II}\bigl((t^{\alpha} - s^{\alpha})^{1/{\alpha}}\bigr).
$$
However the function $\ell$ is not a measure on $\mathcal{B}([0,\infty))$ since for disjoint Borel sets we have here $\ell(A\cup B)^{\alpha} = \ell(A)^{\alpha} + \ell(B)^{\alpha}$.
\vspace{2mm}
In order to give the generating operator for this process at least on the space of polynomials we calculate the following:
\begin{eqnarray*}
\mathbf{E}\bigl( N_{II}(t)^j \big| N_{II}(s) = k\bigr) & = & \sum_{i=0}^j {{j}\choose{i}} k^{j-i} \mathbf{E}\bigl( \bigl(N_{II}(t)- N_{II}(s) \bigr)^i \big| N_s = k\bigr) \\
 & = & k^j  + \sum_{i=1}^j {{j}\choose{i}} k^{j-i} \mathbf{E}\bigl( \bigl(N_{II}(t)- N_{II}(s) \bigr)^i.
 \end{eqnarray*}
 For $i \geqslant 1$ we have
 $$
 \frac{\mathbf{E}\bigl( \bigl(N_{II}(t)- N_{II}(s) \bigr)^i}{t-s} = \frac{t^{\alpha} - s^{\alpha}}{t-s} \sum_{n=1}^{\infty} n^{i-1} \frac{(t^{\alpha} - s^{\alpha})^{n-1}}{(n-1)!} e^{-(t^{\alpha} - s^{\alpha})},
 $$
thus it converges to $\alpha s^{\alpha - 1}$ when $t \rightarrow s$. Consequently we have
$$
\lim_{t\rightarrow s} \frac{\mathbf{E}\bigl( N_{II}(t)^j \big| N_{II}(s) = k\bigr) - k^j}{t-s} = \alpha s^{\alpha - 1} \sum_{i=1}^j {{j}\choose{i}} k^{j-i} = \alpha s^{\alpha - 1} \bigl( (k+1)^j - k^j \bigr).
$$
We see that the generating operator $A_s$ for this process on the polynomial  $W_n(x) = a_0 + a_1 x + \dots + a_n x^n$, $n\in \mathbb{N}$, $a_0, \dots a_n \in \mathbb{R}$ gives the following:
\begin{eqnarray*}
A_s W_n(k) & = &  \lim_{t\rightarrow s} \frac{\mathbf{E}\bigl( W_n\bigl(N_{II}(t)\bigr)\big| N_{II}(s) = k \bigr) - W_n(k)}{t-s} \\
& = & \alpha s^{\alpha - 1} \bigl( W_n(k+1)- W_n(k) \bigr).
\end{eqnarray*}
\subsection{Proper $\vartriangle_{\alpha}$-Generalized Poisson process of type II}
\begin{rem}\label{rem:1}
For the  Kendall convolution the probability kernel for \,$x,y >0$ is given by
$$
h(x,y,t) := \delta_x \vartriangle_{\alpha} \delta_y (0,t) = \left( 1 - \frac{x^{\alpha}y^{\alpha}}{t^{2\alpha}}\right) \mathbf{1}_{\{ x<t,\, y<t\}}.
$$
\end{rem}
\vspace{2mm}
By $F_n$ we denote the cumulative distribution function of $S_n$, $F_1 \equiv F$,   and by $G(t)$ we denote the generalized characteristic function of the step distribution  $\nu$ at the point $t^{-1}$, i.e.
$$
G(t) := h\left( T_{t^{-1}} \nu\right) = \int_0^{\infty} \Bigl( 1 - \frac{x^{\alpha}}{t^{\alpha}} \Bigr)_{+} \nu(dx).
$$

\begin{lem}\label{lem:1}
For the Kendall convolution we have
$$
F_n(t) =  G(t)^{n-1} \bigl[n\left(F(t) - G(t)\right) + G(t) \bigr],
$$
 and
$$
\int_0^t x^{\alpha}dF_n(x) = n t^{\alpha} G(t)^{n-1}\left(F(t) - G(t)\right).
$$
\end{lem}
{\bf Proof.}
Since $G(t)^n = \int_0^{\infty} \bigl( 1 - x^{\alpha}t^{-\alpha} \bigr)_{+} \, dF_n(x)$ then integrating by parts we arrive at
$$
G(t)^n =\alpha t^{-\alpha} \int_0^{t}\! x^{\alpha - 1} F_n(x) dx.
$$
It yields
$$
F_n(t) = \alpha^{-1} t^{1-\alpha} \frac{d}{dt} \bigl[ t^{\alpha} G(t)^n \bigr]= G(t)^n + \frac{nt}{\alpha} G(t)^{n-1} G' (t).
$$
The same procedure applied for $n=1$ leads to the equality $G' (t) = \alpha t^{-1} (F(t) - G(t))$. In order to get the second formula it is enough to notice that
$$
G(t)^n = F_n(t) -  t^{-\alpha} \int_0^{t}\! x^{\alpha} d F_n(x).
$$
\qed

We see now that the distribution of the random variable $N_{II}(t)$ in this case is given by
\begin{eqnarray*}
\mathbf{P} \left\{ N_{II}(t) = 0 \right\} & = &  1 - F(t), \\ \mathbf{P} \left\{ N_{II}(t) = n \right\} & = & G(t)^{n-1} \bigl[ n (F(t) - G(t))(1 - G(t)) + G(t) (1 - F(t)) \bigr].
\end{eqnarray*}
\begin{lem}\label{lem:3}
For the Kendall convolution and for $s<t$
\begin{eqnarray*}
\mathbf{P} \left\{ S_{n+k} <t, S_k <s \right\} = F_n(t) F_k(s) - \frac{s^{\alpha}}{t^{\alpha}} \left(F_n(t) - G(t)^{n}\right) \left( F_k(s) - G(s)^k \right)
\end{eqnarray*}
\end{lem}

\noindent
{\bf Proof.} The result follows from Remark \ref{rem:1} and Lemma \ref{lem:1} by the following calculations
\begin{eqnarray*}
\lefteqn{\mathbf{P} \left\{ S_{n+k} <t, S_k <s \right\} = \int_0^{s}\! \int_0^{t}\! h(x,y,t)  dF_n(x) \, dF_k(y) }\\
 &= &  F_n(t) F_k(s)  - \frac{1}{t^{2\alpha}} \int_0^{t}x^{\alpha} dF_n(x) \int_0^{s}  y^{\alpha}  dF_k(y) \\
 & = & F_n(t) F_k(s) - \frac{s^{\alpha}}{t^{\alpha}} \left(F_n(t) - G(t)^{n}\right) \left( F_k(s) - G(s)^k \right)
\end{eqnarray*}
\qed

\vspace{2mm}

Since the Kendall convolution $\vartriangle_{\alpha}$, $\alpha > 0$, is monotonic and $h(\delta_t) = (1-t^{\alpha})_{+}$ then for the proper  $\vartriangle_{\alpha}$-generalized Poisson process of type II in this case we have $F(t) = \nu([0,t)) = t^{\alpha}\mathbf{1}_{[0,1]}(t) + \mathbf{1}_{(1,\infty)}(t)$. Then
$$
G(t) := \Phi_{\nu}(t^{-1}) = \frac{t^{\alpha}}{2} \mathbf{1}_{[0,1]}(t) + \left( 1 - \frac{1}{2 t^{\alpha}} \right) \mathbf{1}_{(1,\infty)}(t).
$$
Consequently
$$
F_n(t) = \left\{ \begin{array}{lcl}
(n+1) \left(\frac{t^{\alpha}}{2}\right)^n & if & t \leqslant 1; \\[2mm]
\left( 1 + \frac{n-1}{2t^{\alpha}}\right) \left( 1- \frac{1}{2 t^{\alpha}}\right)^{n-1} & if & t > 1.
\end{array}\right.
$$
Now we have
$$
\mathbf{P}\{ N_{II}(t) = 0\} = \left\{ \begin{array}{lcl}
1 - t^{\alpha}  & if & t \leqslant 1; \\
0 & if & t > 1.
\end{array}\right.
$$
and
$$
\mathbf{P}\{ N_{II}(t) = n\} = \left\{ \begin{array}{lcl}
\bigl( n+1 - (n+2) \frac{t^{\alpha}}{2} \bigr) \left( \frac{t^{\alpha}}{2} \right)^n  & if & t \leqslant 1; \\[2mm]
\frac{n}{4t^{2 \alpha}} \left( 1 - \frac{1}{2t^{\alpha}} \right)^{n-1} & if & t > 1. \end{array}\right.
$$
For the convenience we use in the following notation
$$
p:= G(t), \quad q:= G(s).
$$
Then
$$
F_n(t) = \left\{ \begin{array}{ll}
(n+1) p^n & t\leqslant 1, \\
(n(1-p) +p) p^{n-1} & t> 1, \end{array} \right. \quad
\int_0^t x^{\alpha} dF_n(t) = \left\{ \begin{array}{ll}
2n p^{n+1} & t\leqslant 1, \\
\frac{n}{2}  p^{n-1} & t> 1, \end{array} \right.
$$
\begin{eqnarray*}
\lefteqn{\mathbf{P} \left\{ S_{n+k} < t, S_k < s \right\} = } \\[2mm]
 && \hspace{-6mm} \left\{ \begin{array}{ll}
p^{n-1} q^k \bigl[ (n+1)(k+1) p - nkq \bigr]  & s<t<1 \\
p^{n-1} q^{k-1}  \bigl[ nk(1-p)(p-q) + n(1-p)q + k(1-q)p + pq \bigr]  & 1<s<t \\
p^{n-1} q^k  \bigl[ nk(1-p)(1 - 4q(1-p)) + n(1-p) + kp + p\bigr]  &  s<1 < t.
\end{array} \right.
\end{eqnarray*}
In  the case $0<s<t<1$ we have $p=\frac{t^{\alpha}}{2}$, $q=\frac{s^{\alpha}}{2}$,
\begin{eqnarray*}
\lefteqn{\mathbf{P} \left\{ N_{II}(t) = n+k, N_{II}(s) = k  \right\} = } \\[2mm]
 && \left\{ \begin{array}{ll}
 p^{n-2}(p-q) \left[n(p-q)(1-p) + p + q -2p^2 \right] (k+1) q^k& n \geqslant 1, \\[2mm]
 q^k \left[ (k+1)(1-2p +q) -q \right] & n=0,
\end{array} \right.
\end{eqnarray*}
and
\begin{eqnarray*}
\lefteqn{\mathbf{P} \left\{ N_{II}(t)- N_{II}(s) = n \right\} = } \\[2mm]
 &&  \left\{ \begin{array}{ll}
 p^{n-2} \frac{p-q}{(1-q)^2} \left[ n(p-q)(1-p)+ p + q -2p^2 \right] & n \geqslant 1, \\[2mm]
 1 - \frac{2(p-q)}{(1-q)^2} & n=0.
 \end{array} \right.
\end{eqnarray*}
In particular for $n\geqslant 1$ we have that
\begin{eqnarray*}
\lefteqn{\mathbf{P} \left\{ N_{II}(t) - N_{II}(s) = n, N_{II}(s) = k  \right\}} \\
 & = & \mathbf{P} \left\{ N_{II}(t)- N_{II}(s) = n \right\} \, (k+1) (1-q)^2 q^k.
\end{eqnarray*}
Since in this case $\mathbf{P} \left\{ N_{II}(s) = k  \right\} = ( k (1-q) +1-2q)q^k$ we see that the increments of this process are not independent. The dependence of the distribution of increment $N_{II}(t) - N_{II}(s)$ on the time increment $[s,t)$ is also rather complicated and far from linear. \\

Notice also that the process $\{ N_{II}(t) \colon t \geqslant 0\}$ is not a Markov process. To see this it is enough to check the simplest version of Markov condition: for $1>t>s>u$, $k \in \mathbb{N}_0$ with the notation $r = {{u^{\alpha}}/2}$ we have
\begin{eqnarray*}
\lefteqn{\mathbf{P} \left\{ N_{II}(t) = k \big| N_{II}(s) = k, N_{II}(u) = k \right\}} \\
 & = & \frac{\mathbf{P} \{ S_k <u \} - \mathbf{P} \{ S_k < u, S_{k+1} <t \}}{\mathbf{P} \{ S_k <u \} - \mathbf{P} \{ S_k < u, S_{k+1} < s\}} = \frac{(k+1)(1-2p) +kr}{(k+1)(1-2q) +kr}.
\end{eqnarray*}
On the other hand
\begin{eqnarray*}
\lefteqn{\mathbf{P} \left\{ N_{II}(t) = k \big| N_{II}(s) = k \right\}} \\
 & = & \frac{\mathbf{P} \{ S_k <s \} - \mathbf{P} \{ S_k < s, S_{k+1} <t \}}{\mathbf{P} \{ S_k <s \} - \mathbf{P} \{ S_{k+1} < s\}} = \frac{(k+1)(1-2p) +kq}{(k+1)(1-2q) +kq}.
\end{eqnarray*}
For the completeness we calculated also the remaining part of  distribution, thus in the case $1<s<t$ we have $p=1 - \frac{1}{2t^{\alpha}}$, $q=1 - \frac{1}{2s^{\alpha}}$,
\begin{eqnarray*}
\lefteqn{\mathbf{P} \left\{ N_{II}(t) = n+k, N_{II}(s) = k  \right\} = } \\[2mm]
 &&  \left\{ \begin{array}{ll}
 p^{n-2}(p-q)(1-p)^2 \left[ n(p-q) + p + q \right] k q^{k -1} & n \geqslant 1, \\[2mm]
 (1-p)^2 k q^{k-1} & n=0,
 \end{array} \right.
\end{eqnarray*}
and
\begin{eqnarray*}
\lefteqn{\mathbf{P} \left\{ N_{II}(t)- N_{II}(s) = n \right\} = } \\[2mm]
 &&  \left\{ \begin{array}{ll}
 p^{n-2} (p-q) \frac{(1-p)^2}{(1-q)^2}  \left[ n(p-q) + p + q \right] & n \geqslant 1, \\[2mm]
 \frac{(1-p)^2}{(1-q)^2} & n=0.
 \end{array} \right.
\end{eqnarray*}
In the case $s<1<t$ we have $p=1 - \frac{1}{2t^{\alpha}}$, $q=\frac{s^{\alpha}}{2}$,
\begin{eqnarray*}
\lefteqn{\hspace{-3mm}\mathbf{P} \left\{ N_{II}(t) = n+k, N_{II}(s) = k  \right\} = p^{n-2} q^k(1-p)^2 } \\[2mm]
 &&  \hspace{-6mm}\left\{ \begin{array}{ll}
 \bigl[n k (p-q)(1 - 4q(1-p)) + n (1-2q) (p-2q(1-p)) \\[2mm]
 \hspace{3mm} + kq(1 - 4q +4p^2) + 2q (1-2q) \bigr]  & n \geqslant 1, \\[2mm]
 4 k  p^2 q  & n=0,
 \end{array} \right.
\end{eqnarray*}
and
\begin{eqnarray*}
\lefteqn{\mathbf{P} \left\{ N_{II}(t)- N_{II}(s) = n \right\} =  p^{n-2} \frac{(1-p)^2}{(1-q)^2}\, \times} \\[2mm]
 && \hspace{-3mm}\left\{ \begin{array}{ll}
    \bigl[ n\bigl(4q^2 (1-p)^2 +(1-q)^2 - (1-p)\bigr) + q \bigl( 4qp^2 + 2 - 5q\bigr) \bigr] & n \geqslant 1, \\[2mm]
 4p^2 q^2  & n=0.
 \end{array} \right.
\end{eqnarray*}
Now we are able to calculate the conditional distributions and conditional moments. Notice first that
$$
\sum_{n=0}^{\infty} \mathbf{P} \left\{ N_{II}(t) = n+k, N_{II}(s) = k  \right\} = \mathbf{P} \left\{ N_{II}(s) = k  \right\}.
$$
Consequently
\begin{eqnarray*}
\lefteqn{\mathbf{E}\left( N_{II}(t)^j \big| N_{II}(s) = k \right) - k^j } \\
& = & \sum_{n=1}^{\infty} \left( (n+k)^j - k^j \right)\mathbf{P} \left\{ N_{II}(t) = n+k\big| N_{II}(s) = k  \right\}.
\end{eqnarray*}
Calculating this value we will need the following sequence of functions
$$
\Phi_j (k, p) = \sum_{n=1}^{\infty} (n+k)^j p^{n+k},
$$
where $k \in \mathbb{N}$ and the natural domain for $p$ is $(-1,1)$, however we need only $p\in (0,1)$. It is easy to see that
$$
\Phi_{j+1} (k, p) = p \, \frac{\partial}{\partial p} \Phi_j (k, p),
$$
which is a pretty simple recursive formula, however it leads to rather laborious calculations.
In particular we have
\begin{eqnarray*}
&& \Phi_0 (k, p)= \frac{p^{k+1}}{(1-p)}, \quad \Phi_1 (k, p) = \frac{k(1-p) +1}{(1-p)^2} \, p^{k+1}, \\[2mm]
&& \Phi_2 (k, p) = \frac{k^2(1-p)^2 + 2k(1-p) + 1 + p}{(1-p)^3} \, p^{k+1}.
\end{eqnarray*}
Notice also that we have
$$
\sum_{n=1}^{\infty} n (n+k)^j p^{n+k} = \Phi_{j+1} (k, p) - k\Phi_j (k, p).
$$
{\bf 1) In the case $0<s<t<1$} we have
\begin{eqnarray*}
\lefteqn{\mathbf{E}\left( N_{II}(t)^j \big| N_{II}(s) = k \right) - k^j = \frac{(k+1)}{(k(1-q)+1-2q)} \, \frac{p-q}{p^{k+2}} \times }\\
 &&\hspace{-8mm} \Bigl[ (1-p)(p-q) \bigl( \Phi_{j+1}(p) - k \Phi_j(p)\bigr) + (p+q - 2p^2) \Phi_j(p)
   - 2 k^j p^{k+2} \Bigr].
\end{eqnarray*}
Since $\lim_{t\rightarrow s} \frac{p-q}{t-s} = \frac{1}{2} (s^{\alpha})' = \frac{1}{2} \alpha s^{\alpha -1}$ we finally obtain
\begin{eqnarray*}
\lefteqn{\lim_{t \searrow s} \frac{ \mathbf{E}\left( N_{II}(t)^j \big| N_{II}(s) = k \right) - k^j}{t-s} }\\
& = & \frac{(k+1)\alpha s^{\alpha - 1} q^{-1-k}}{(k+1)(1-q) -q}  \Bigl[ (1-q) \Phi_j(q)  - k^j q^{k+1}  \Bigr].
\end{eqnarray*}
Notice that if $ s \nearrow 1$ then this value converges to
$$
\frac{k+1}{k} \, \alpha 2^{k+1} \left[ \Phi_j\Bigl(k, \frac{1}{2}\Bigr) - k^j 2^{-k} \right].
$$
{\bf 2) In the case $1<s<t$} we have
\begin{eqnarray*}
\lefteqn{\mathbf{E}\left( N_{II}(t)^j \big| N_{II}(s) = k \right) - k^j = \frac{(1-p)^2}{(1-q)^2} \, \frac{(p-q)}{p^{k+2}} \times }\\
 && \hspace{-8mm}\left[ (p-q) \Bigl( \Phi_{j+1}(k,p) - k \Phi_j(k,p)\bigr) + (p+q) \Phi_{j}(k,p) - \frac{k^j p^{k+2}}{(1-p)^2} (2-p-q)\right].
\end{eqnarray*}
Notice that $\lim_{t\rightarrow s} \frac{p-q}{t-s} = \frac{\alpha}{2 s^{\alpha+1}} $ where $q = 1 - \frac{1}{2s^{\alpha}}$. Consequently
$$
\lim_{t\searrow s} \frac{\mathbf{E}\left( N_{II}(t)^j \big| N_{II}(s) = k \right) - k^j}{t-s} = \frac{\alpha s^{\alpha -1}}{s^{2\alpha} q^{k+1}} \Bigl( \Phi_{j}(k,q) -  \frac{ k^{j} q^{k+1} }{1-q} \Bigr).
$$
For $s \searrow 1$ this value converges to
$$
\alpha 2^{k+1} \Bigl( \Phi_{j}\Bigl(k,\frac{1}{2}\bigr) -  k^{j}2^{-k } \Bigr)
$$
which is different from the analogous limit from the left hand side of $1$.

\vspace{2mm}

{\bf 3) In the case $s<1<t$} we have
\begin{eqnarray*}
\lefteqn{\mathbf{E}\left( N_{II}(t)^j \big| N_{II}(s) = k \right) - k^j = \frac{(1-p)^2 p^{-k-2}}{k(1-q) +1 - 2q}\, \times }\\
 && \hspace{-7mm}\Bigl[ \big[ k (1-4q + 4p^2) +2(1-2q) \big] q \Phi_j \\
 && \hspace{-7mm} + \bigl[(1-2q)(p-2q(1-p)) +k(p-q)(1-4q(1-p)) \bigr]\bigl( \Phi_{j+1} - k\Phi_j\bigr) \\
 && \hspace{-7mm} - \, \frac{k^j p^{k+2}}{(1-p)^2}\bigl[1-2q + k\bigl( 1 -q - 4q (1-p)^2 \bigr) \bigr] \Bigr],
 \end{eqnarray*}
 where $\Phi_j := \Phi_j(k,p)$.

For every $s>0$, $s\neq 1$ let $\Gamma_s$ be a random variable with the distribution
$$
\mathbf{P} \left\{\Gamma_s = n \right\} = (1-q) q^{n-1}, \quad n=1,2,\dots, \quad q = G(s).
$$
We see that in spite of the fact that the process $\{ N_{II}(t)\colon t \geqslant 0\}$ does not have Markov property still the analogue of infinitesimal operator $A_s$, $s\neq 1$ is well defined on the set
$$
\left\{ f \colon \mathbb{N} \rightarrow \mathbb{R}: \mathbf{E}|f(\Gamma_s +k)| < \infty \,\, \forall\, k \in \mathbb{N}_0 \right\}
$$
and
\begin{eqnarray*}
A_s f(k) & \stackrel{def}{=} &  \lim_{t\rightarrow s} \frac{\mathbf{E}\bigl( f\bigl(N_{II}(t)\bigr)\big| N_{II}(s) = k \bigr) - f(k)}{t-s} \\
& = & C(k,s) \bigl( \mathbf{E}f(\Gamma_s +k) - f(k) \bigr),
\end{eqnarray*}
where
$$
C(k,s) = \left\{ \begin{array}{ll}
\frac{2(k+1) \alpha s^{\alpha - 1}}{(k+1)(2-s^{\alpha}) - s^{\alpha}} & s<1; \\[2mm]
2\alpha s^{-1} & s>1.
\end{array}\right.
$$

\end{document}